\documentclass[12pt,leqno]{article}

\usepackage{amsfonts}
\usepackage{amsthm}
\usepackage[mathscr]{eucal}
\usepackage{amsmath}
\usepackage{amssymb}
\usepackage{amscd}

\makeatletter
\def\newtheorem#1{\@ifnextchar[{\@othm{#1}}{\@nthm{#1}}}
\def\@nthm#1#2{%
\@ifnextchar[{\@xnthm{#1}{#2}}{\@ynthm{#1}{#2}}}
\def\@xnthm#1#2[#3]{\expandafter\@ifdefinable\csname #1\endcsname
{\@definecounter{#1}\@newctr{#1}[#3]%
\expandafter\xdef\csname the#1\endcsname{\expandafter\noexpand
  \csname the#3\endcsname \@thmcountersep \@thmcounter{#1}}%
\global\@namedef{#1}{%
  \@thm{#1}{#2}}\global\@namedef{end#1}{\@endtheorem}}}
\def\@ynthm#1#2{\expandafter\@ifdefinable\csname #1\endcsname
{\@definecounter{#1}%
\expandafter\xdef\csname the#1\endcsname{\@thmcounter{#1}}%
\global\@namedef{#1}{%
  \@thm{#1}{#2}}\global\@namedef{end#1}{\@endtheorem}}}
\def\@othm#1[#2]#3{%
  \@ifundefined{c@#2}{\@nocounterr{#2}}%
  {\expandafter\@ifdefinable\csname #1\endcsname
  {\global\@namedef{the#1}{\@nameuse{the#2}}%
\global\@namedef{#1}{\@thm{#2}{#3}}%
\global\@namedef{end#1}{\@endtheorem}}}}
\def\@thm#1#2{\refstepcounter
    {#1}\@ifnextchar[{\@ythm{#1}{#2}}{\@xthm{#1}{#2}}}
\def\@xthm#1#2{\@begintheorem{#2}{\csname the#1\endcsname}\ignorespaces}
\def\@ythm#1#2[#3]{\@opargbegintheorem{#2}{\csname
       the#1\endcsname}{#3}\ignorespaces}
\def\@thmcounter#1{\noexpand\arabic{#1}}
\def\@thmcountersep{.}
\def\@begintheorem#1#2{\trivlist
   \item[\hskip \labelsep{\bfseries #2\ #1.}]\itshape}
\def\@opargbegintheorem#1#2#3{\trivlist
      \item[\hskip \labelsep{\bfseries #2\ #1\ [#3].}]\itshape}
\def\@endtheorem{\endtrivlist}
\makeatother
\theoremstyle{plain}
\newtheorem{lemma}{Lemma}[section]
\newtheorem{theorem}[lemma]{Theorem}

\theoremstyle{definition}

\newtheorem{definition}[lemma]{Definition}

\def\d{$\displaystyle}
\def\be{\begin{equation}}
\def\heavycdot{\raisebox{2pt}{\tiny\kern.5em$\bullet$}}

\numberwithin{equation}{section}
\def\begeq{\stepcounter{lemma}\begin{equation}}

\date{}

\begin{document}

\title{Homology of $L_\infty$-algebras and Cyclic Homology}
\author{Masoud Khalkhali\thanks{\; The author is supported by NSERC of
Canada.  The author is grateful to J. Williams for typing the
manuscript.}}

\maketitle



A celebrated theorem of Loday and Quillen [LQ] and (independently)
Tsygan [T] states that the Lie algebra homology of the Lie algebra of
stable matrices over an associative algebra is canonically isomorphic,
as a Hopf algebra, to
the exterior power of the cyclic homology of the associative algebra.
The main point of this paper is to lay the ground such that an extension
of this theorem to the category of $A_\infty$-algebras becomes possible
(theorem 3.1).

The category of $L_\infty$- (respectively, $A_\infty$-) algebras extend
the category of differential graded $(DG)$ Lie (respectively, $DG$
associative) algebras.  These concepts are both due to J. Stasheff.  See
[S], [LS], and references therein, and also [HS] where an alternative
approach to $L_\infty$-algebras is given.  In [Kh], we proposed an
approach to homological invariants of $A_\infty$-algebras (Hochschild,
cyclic, periodic cyclic, etc.) based on the notion of $X$-complex due to
Cuntz and Quillen [CQ].  It seems that it is now possible to extend most
of the tools of noncommutative geometry of Connes [C] to the homotopical
setting of $A_\infty$ and $L_\infty$-algebras.  There is, however, a
notable exception in that so far we don't know how the $K$-theory of an
$A_\infty$-algebra should be defined.

\section{$L_\infty$ and $A_\infty$ algebras}
Let $V$ be a vector space (not graded).  Let $S^cV$ denote the cofree
cocommutative counital coassociative coalgebra generated by $V$.  Over
fields of characteristic zero there are two different constructions for
$S^cV$ that we recall now.  First, let \d T^cV = \bigoplus_{n\geq
0}V^{\otimes n}$ be the cofree counital coassociative coalgebra
generated by $V$.  Its comultiplication $\Delta :T^cV\longrightarrow
T^cV\otimes T^cV$ is defined by
\[ \Delta (v_1,\dots ,v_n) = \sum^n_{i=0}(v_1,\dots ,v_i)\otimes
(v_{i+1},\dots ,v_n)\; , \]
and its counit is the projection onto $V^{\otimes0}=K$.  The
symmetric group $S_n$ acts (on the left) on $V^{\otimes n}$ in the
standard way (without signs):
\[ \sigma (v_1,\dots ,v_n) = (v_{\sigma^{-1}_{(1)}},\dots
,v_{\sigma^{-1}_{(n)}})\; . \]
Let
\[ S^cV=\bigoplus_{n\geq 0}(V^{\otimes n})^{S_n} \]
be the space of invariants of this action (symmetric tensors).  Since
the $S_n$- action is compatible with the coproduct, $S^cV$ is a
subcoalgebra of $T^cV$.  It is obviously cocommutative.  The
symmetrization map
\[ P(v_1,\dots ,v_n)=\frac{1}{n!}\sum_{\sigma\in S_n}\sigma (v_1,\dots
,v_n)\; , \]
defines a surjection $T^cV\to S^cV$.  It is not a coalgebra map, but is
a retraction for the canonical inclusion $S^cV\to T^cV$.

Alternatively, let \d \overline{S}^cV=\bigoplus_{n\geq 0}(V^{\otimes
n})_{S_n}$ be the space of coinvariants of the above $S_n$-action.  We
have a canonical isomorphism of vector spaces
$S^cV\simeq\overline{S}^cV$.  In fact the natural map $(V^{\otimes
n})^{S_n}\to (V^{\otimes n})_{S_n}$ obtained by composing the inclusion
into $V^{\otimes n}$ by projecting onto coinvariants has an inverse
given by the symmetrization $P$.  Under this isomorphism the coproduct
of $S^cV$ is transformed to a coproduct $\overline{\Delta}$ an
$\overline{S}^cV$.  It is easy to see that it is given by
\[ \overline{\Delta}(v_1,\dots ,v_n)=\sum_{\sigma\in
S_{p,n-p}}\sigma\Delta (v_1,\dots ,v_n), \]
where $S_{p,q}$ is the space of all $(p,q)$-shuffles.

The coalgebra $S^cV$ has universal properties with respect to coalgebra
morphisms and coderivations.  Let $A$ be a cocommutative counital
coalgebra.  Then we have a 1-1 correspondence
\[ Hom^{coalg}(A,S^cV)\simeq Hom(A,V)\; , \]
where $Hom^{coalg}$ denotes counital coalgebra maps.

Let $M$ be a left counitary comodule over a cocommutative coalgebra $A$,
with left comultiplication $\Delta_\ell :M\to A\otimes M$.  Let
$\Delta_r:M\to M\otimes A$ be the natural right coaction derived from
$\Delta_\ell$.  Recall that a coderivation is a linear map $D:M\to A$
such that
\[ \Delta D=(1\otimes D)\Delta_\ell +(D\otimes 1)\Delta_r\; . \]
For any (symmetric) $S^cV$ comodule $M$, we have a 1-1 correspondence
\[ coder(M,S^cV)\simeq Hom(M,V)\; . \]

A coderivation $\delta :T^cV\to T^cV$ induces a coderivation
$\delta^\prime :S^cV\to S^cV$ by the composition of maps
\[ S^cV\longrightarrow T^cV\overset{\delta}{\longrightarrow}
T^cV\longrightarrow S^cV \]
where the first map is the inclusion and the last map is the projection
$P$.

\begin{lemma}
$\delta^\prime$ is a coderivation.
\end{lemma}

Instead of a formal proof, it is
best to observe that this is the (pre) dual of a fact for algebras:
a derivation $\delta :TV\to TV$ of free algebras induces a derivation
$\delta^\prime :SV\to SV$ of the symmetric algebra via the composition
\[ SV\longrightarrow TV\overset{\delta}{\longrightarrow}TV\longrightarrow
SV\; . \]
In fact the canonical projection
$TV\longrightarrow SV$ has a natural section
$SV\overset{s}{\longrightarrow}TV$, the symmetrization map.
$s$ is not an algebra map, but $s(xy)-s(x)s(y)\in I$ the ideal generated
by $\{ u\otimes v-v\otimes u;u,v\in V\}$ in $TV$.  Hence to show that
$\delta^\prime$ is a derivation, suffices to show that $\delta I\subset
I$.  However, for an element $z = x(u\otimes v-v\otimes u)y\in I$, we
have
\[ \delta z=\delta x\cdot (u\otimes v-v\otimes u)y+x(\delta u\otimes
v+u\otimes\delta v-\delta v\otimes u-v\otimes\delta u)y+x(u\otimes
v-v\otimes u)\delta y \]
and hence $\delta z\in I$.

Next, we note that all of our definitions and constructions so far are
functorial and can be repeated verbatim in any symmetric monoidal
category, also known as symmetric tensor category.  In particular we can
apply this to the tensor category of $\mathbb{Z}$-graded vector space.

The objects are graded vector spaces \d V =
\bigoplus_{i\in\mathbb{Z}}V_i$ and the morphisms \d Hom(V,W) =
\bigoplus_{n\in\mathbb{Z}}Hom^n(V,W)$, where \d Hom^n(V,W) =
\prod_{i\in\mathbb{Z}}Hom(V_i,W_{i+n})$ is the set of linear maps of
degree n.  The tensor product is defined by \d (V\otimes W)_n =
\bigoplus_{i+j=n}V_i\otimes W_j$.  The symmetry $S:V\otimes W\to
W\otimes V$ is defined by $S(v\otimes w)=(-1)^{|v|\; |w|}w\otimes v$.
As in any symmetric tensor category, the symmetric group $S_n$ acts on
the tensor power $V^{\otimes n}$.  It is given by
\[ \sigma (v_1,\dots ,v_n)=(-1)^\varepsilon (x_{\sigma^{-1}(1)},\dots
,x_{\sigma^{-1}(n)})\; , \]
where $(-1)^\varepsilon$ is, in general, different from the sign of the
permutation $\sigma$.  It is however the same if $|v_i|=1$ for all $i$.
>From now on we freely use our previous constructions in the above
context.

In particular, we note that the cofree graded cocommutative coalgebra of
a vector space $V$ (not graded), denoted $\Lambda^cV$, is nothing but
$S^cV[1]$.  Motivated by this, we define for any positively graded
vector space \d V=\bigoplus^\infty_{i=0}V_i$, the cofree cocommutative
coalgebra $\Lambda^cV:=S^cV[1]$.  We have $\Lambda^0V=K$,
$\Lambda^nV=\{x\in V[1]^{\otimes m};\sigma x=x\}$.  The coproduct is given by
\[ \Delta (v_1,\dots ,v_n) = \sum^n_{i=1}(v_1,\dots ,v_i)\otimes
(v_{i+1},\dots ,v_n)\; . \]

Under the graded bracket $[\delta^1,\delta^2] =
\delta^1\delta^2-(-1)^{|\delta^1|\; |\delta^2|}\delta^2\delta^1$, the
space $coder(\Lambda^cV,\Lambda^cV)$ is a graded Lie algebra.  The
isomrophism
\[ coder(\Lambda^cV,\Lambda^cV) = Hom(\Lambda^cV,V[1]) \]
defines a bracket on the right hand side.  An element $\ell
:\Lambda^cV\to V[1]$ of degree -1 is naturally identified with a
sequence of map
\[ \ell_n:V[1]^{\otimes n}\to V[1],\qquad n\geq 1, \]
such that $|\ell_n|=-1$ and $\ell_n$ is ``antisymmetric'', in the sense
that, $\forall\sigma\in S_n$,
\[ \ell_n(v_{\sigma (1)},\dots ,v_{\sigma (n)}) = (-1)^\varepsilon
\ell_n(v_1,\dots ,v_n), \]
where $\varepsilon$ is as before, except that we now use degrees in
$V[1]$.

Given such a map $\ell$, the corresponding coderivation $\delta_\ell
:\Lambda^cV\to\Lambda^cV$ is given by the
\[ \delta_\ell (v_1,\dots ,v_n) = \sum^n_{k=1}\;\sum^{n-k}_{j=1}
(-1)^\varepsilon P(v_1,\dots ,\ell_k(v_j,\dots ,v_{j+k}),\dots ,v_n)\; , \]
where $P$ is the projection operator and \d \varepsilon =
\sum^{j-1}_{s=1}|v_s|$ (suspended degrees).

\begin{definition}
Let \d L = \bigoplus^\infty_{i=0}L_i$ be a graded vector space.  An
$L_\infty$ structure (also called strongly homotopy Lie algebra
structure) on $L$ is a coderivation $\delta_\delta
:\Lambda^cL\to\Lambda^cL$ of degree -1 such that $\delta^2_\ell =0$.
\end{definition}

Since $\delta^2_\ell =\frac12 [\delta_\ell ,\delta_\ell ] =
\frac12\delta_{[\ell\cdot\ell ]} = \delta_{\ell\circ\ell}$.  We see that
$\delta^2_\ell=0$ iff $\ell\circ\ell =0$.  Here $\circ$ is the analogue
of Gerstenhaber's product for Lie cochains.

One can give a completely parallel treatment for $A_\infty$-algebras.
>From the beginning one has to work with cofree coassociative coalgebra
$T^cV$.  Thus an $A_\infty$ structure on a graded vector space \d
A=\bigoplus_{i\geq 0}A_i$ is a degree -1 coderivation
$\delta_m:T^cA[1]\to T^cA[1]$ such that
$\delta^2_m=\frac12[\delta_m,\delta_m]=0$.  Here \d
m=\sum^\infty_{i=1}m_i:T^cV[1]\to V[1]$ is a cochain of degree -1.  Each
individual cochain $m_i:V[1]^{\otimes^i}\to V[1]$ has degree -1.

The functor $A\to A^{Lie}$ from associative algebras to Lie algebras,
where $A^{Lie}=A$ with the Lie bracket $[a,b]=ab-ba$, has a vast
generalization to a functor from $A_\infty$-algebras to
$L_\infty$-algebras that we are going to describe now.  In fact this is
just the analogue of lemma 1.1 above in the tensor category of graded
vector spaces.  Note that a direct computation
would be extremely hard.  Let $(A,m)$ be an $A_\infty$-algebra.  From
the coderivation $\delta_m:T^cA[1]\to T^cA[1]$, we construct a coderivation
$\delta_\ell :\Lambda^cA=S^cA[1]\to\Lambda^cA$ as the composition
\[ S^cA[1]\overset{i}{\hookrightarrow}
T^cA[1]\overset{\delta_\ell}{\longrightarrow}
T^cA[1]\overset{p}{\longrightarrow} S^cA[1] \]
where the first map is the natural inclusion and the last map is the
natural projection so that $\delta_\ell =p\delta_mi$.  Since $ip=id$, we
have
\[ \delta^2_\ell =p\delta_mip\delta_mi=p\delta^2_mi=0\; . \]

We also need to define the tensor product of two $A_\infty$-algebras
$(A,m)$ and $(A^\prime ,m)$.  This is an $A_\infty$-algebra whose
underlying graded space is $A\otimes A^\prime$.  Our construction of
this tensor product is based on the following simple facts.  First, let
$C_1$ and $C_2$ be graded coalgebras and $\delta_i:C_i\to C_i,\; i=1,2$
graded coderivations of degree -1 such that $\delta^2_1=\delta^2_2=0$.
Then the map
\[ \delta=\delta_1\otimes 1+1\otimes\delta_2: C_1\otimes C_2\to C_1\otimes
C_2 \]
is a coderivation and $\delta^2=0$.  Note that graded tensor product of
graded coalgebras is used.

Secondly, if $f:C_1\to C_2$ is a coalgebra map and $\delta :C_2\to C_2$
is a coderivation, then $\delta\circ f:C_1\to C_2$ is a coderivation,
where we regard $C_1$ a $C_2$-bicomodule via $f$.  There is a similar
statement for composing from right, i.e. $f\circ\delta$.

Now let $(A,m)$ and $(A^\prime ,m)$ be $A_\infty$-algebras.  Consider
the composition of maps
\[ \begin{CD} T^c(A\otimes A^\prime )[1]\longrightarrow T^cA[1]\otimes
T^cA^\prime[1] @>{\delta_m\otimes
1+1\otimes\delta_m}>> T^cA[1]\otimes
T^cA[1]\longrightarrow T^c(A\otimes A^\prime )[1] \end{CD} \]
where the first and last maps are morphisms of graded coalgebras.  The
first map is defined by a pair of cocommuting coalgebra maps
$T^c(A\otimes A^\prime )[1]\to T^cA[1]$ and $T^c(A\otimes A^\prime )[1]\to
T^cA^\prime [1]$,
while each of these individual maps are defined by sending $(A\otimes
A^\prime )[1]\to A[1]$ (resp. $(A\otimes A^\prime )[1]\to A^\prime [1]$)
and other parts of $T^c(A\otimes A^\prime )[1]$ to zero and then
coextending them to a coalgebra map.  The last map is defined in a
similar way.  Of course one has to check that the induced $T^c(A\otimes
A^\prime )[1]$-bicomodule structure on itself, induced from these two
maps is the same as the bicomodule-structure induced from the coproduct.

There is a special case of the above construction that is particularly
simple and also important for us.  If $(A,m)$ is an $A_\infty$-algebra
and $B$ is an associative algebra considered as an $A_\infty$-algebra in
the obvious way, then $A\otimes B$ is given by $(A\otimes B)_n =
A_n\otimes B$ and its defining cochains
$m^\prime_n:(A\otimes B)^{\otimes n}\to A\otimes B$ are given by
\[ m^\prime_n(a_1\otimes b_1,\dots ,a_n\otimes b_n) =
m_n(a_1,\dots ,a_n)\otimes b_1b_2,\dots ,b_n\; . \]
In particular, we will apply this construction to $B=M_n(\mathbb{K})$, the
algebra of $n\times n$ matrices over the ground field $\mathbb{K}$.  We
obtain an $A_\infty$-algebra $M_n(A)$.  We denote its associated
$L_\infty$-algebra by $g\ell_n(A)$, and we let $g\ell
(A)=\underset{\longrightarrow}{\lim}\; g\ell_n(A)$ be the direct limit of
these $L_\infty$-algebras.

\section{Homology of $L_\infty$ and $A_\infty$ algebras}
Like Lie algebras, there are at least two approaches to define the
(co)homology of an $L_\infty$-algebra $(L,\ell )$ with coefficients in
an $L_\infty$-module.  One can define a universal enveloping algebra
$U(L,\ell )$ and define the (co)homology of $(L,\ell )$ as the Hochschild
(co)homology of $U(L,\ell)$ with coefficients.  Note that there are two
variants of this universal enveloping algebra, one of which is an
$A_\infty$-algebra and can be regraded as the left adjoint to the
functor $A\to A^{Lie}$ from $A_\infty$ to $L_\infty$-algebras.  The
second is a DG algebra [LS, HS].  One obviously hopes that the two definitions
give the same answer.

Alternatively, one notes that $(\Lambda^cL,\delta_\ell )$ is the analogue
of the Chevally-Eilenberg complex for $L_\infty$ algebras.  In fact if
$L$ is a (DG) Lie algebra, then it is exactly the Chevally-Eilenberg
complex of $L$.  Again one naturally expects that this definition be
equivalent to the above.  We, however, don't need this and indeed are
only interested in homology with trivial coefficients.  We denote this
by $H_\bullet (L,\ell )$ and define it as the homology of graded
cocommutative coalgebra $(\Lambda^cL,\delta_\ell )$.  In
particular $H_\bullet (L,\ell )$ is a graded cocommutative coalgebra.

We need to know that ``inner derivations'' act like zero on $H_\bullet$.
Let $(L,\ell )$ be an $L_\infty$-algebra.  By a derivation of $L$ of
degree $k$, we mean a map
\[ d:\Lambda^cL\to L[1] \]
such that $|d|=k$ and the induced coderivation
\[ \delta_d:\Lambda^cL\to\Lambda^cL \]
satisfies $[\delta_\ell ,\delta_d]=0$.  Thus a derivation induces a
morphism of the complex $(\Lambda^cL,\delta_\ell )$ and hence act on
$H_\bullet (L,\ell )$.  To give an example, let $L$ be a Lie algebra.
A map $d_k:\Lambda^kL\to L$ and $d_i=0$, $i\not= k$, is a derivation in
the above sense iff $L$ is a Chevalley-Eilenberg cocycle for $L$ with
coefficients in the adjoint module $L$.

A derivation $\delta_d$ of $(L,\ell )$ is called {\em inner} if there
exists a cochain $\delta_{d^\prime}$ such that $\delta_d=[\delta_\ell
,\delta_{d^\prime}]$.

\begin{lemma}
Inner derivations act like zero on the homology of $L_\infty$-algebras.
\end{lemma}

Now let $(L,\ell )$ be an $A_\infty$-algebra and $h\subset L_0$ a sub Lie
algebra.  One can define a reduced Chevalley-Eilenberg complex,
$CE_\bullet (L,\ell )_h$.  It is simply the space of coinvariants
$CE_\bullet (L,\ell )$ under the action of $h$ by inner derivations.  It
follows from the above lemma, exactly as in the Lie algebra case (cf.
$[L]$), that if $h$ is a reductive Lie algebra, then the projection
$CE_\bullet (L,\ell )\to CE_\bullet (L,\ell )_h$ is a
quasi-isomormphism.

Next, let us briefly recall the definition of the cyclic homology of
$A_\infty$-algebras.  The cyclic homology of $A_\infty$-algebras was
first defined by Getzler and Jones in [GJ], where they defined a bicomplex
similar to Connes's $(b,B)$-bicomplex.

In [Kh], we gave an alternative approach to cyclic homology of
$A_\infty$-algebras, based on the notion of $X$-complex due to Cuntz and
Quillen [CQ].  In particular a new bicomplex for $A_\infty$-algebras is
defined in [Kh], which is the analogue of Connes-Tsygan bicomplex.  For this
paper, however, it is most convenient to use a third complex, the analogue
of Connes's cyclic complex.  It follows from the results of [Kh], that
the three complexes are quasi-isomorphic.

Let $(A,m)$ be an $A_\infty$-algebra.  For $n\geq 0$, let $C_n(A,m) =
A[1]^{\otimes n+1}$.  Define operators $b_m:C_n(A,m)\to C_{n-1}(A,m)$ and
$\lambda :C_n(A,m)\to C_n(A,m)$ by
\begin{eqnarray*}
b_m(a_1,\dots ,a_n) & = &
\sum^n_{i=1}\;\sum^{n-i}_{j=1}(-1)^{\varepsilon_{ij}}(a_1,\dots
,m_i(a_j,\dots ,a_{j+1}),\dots ,a_n)\\
& + & \sum^n_{i=1}\;\sum^n_{j=n-i+1}(-1)^\varepsilon
(m_i(a_j,\dots ,a_{i+j-n-2}),\dots ,a_{j-1})
\end{eqnarray*}
and
\[ \lambda (a_1,\dots ,a_n) = (-1)^{|a_n|(|a_1|+\cdots +|a_{n-1}|)}
(a_n,a_1,\dots ,a_{n-1})\; . \]
Then $b^2_m=0$.  Also $Im(1-\lambda )$ is invariant under $b_m$.  Let
$C^\lambda_n(A,m)=C_n(A,m)/Im(1-\lambda )$ be the space of coinvariants
under the action of the cyclic group.  The complex $(C^\lambda_\bullet
(A,m),b_m)$ is called the cyclic complex of the $A_\infty$ algebra
$(A,m)$ and its homology is called the cyclic homology of $(A,m)$ and is
denoted by $HC_\bullet (A,m)$.

\section{Proof of the main theorem}
Let $(A,m)$ be a unital $A_\infty$-algebra and let $(g\ell (A),\ell )$
denote the $L_\infty$-algebra of stable matrices over $(A,m)$.  Let
$g\ell (K)\subset g\ell (A)$ be the subalgebra of stable matrices with
coefficients in the ground field $K$.  Our first task is to show that
the complex of coinvariants $CE_\bullet (g\ell (A),\ell )_{g\ell (K)}$ is a
cocommutative and commutative DG Hopf algebra, as in the associative
case.

In fact, it is enough to observe that the same maps as in [L] serves the
purpose.  Define
\[ \oplus :g\ell (A)\oplus g\ell (A)\to g\ell (A) \]
by sending $(a_{ij})\times (b_{ij})\mapsto (c_{ij})$, where
\[ c_{2i+1,2j+1}=a_{ij},\quad c_{2i,2j}=b_{ij}\quad\mbox{and}\quad
c_{ij}=0\quad\mbox{otherwise.} \]
It is easy to see that $\oplus$ is a strict morphism of
$L_\infty$-algebras and hence, via the diagonal map, induces a morphism
of complexes
\[ CE_\bullet (g\ell (A),\ell )\otimes CE_\bullet (g\ell (A),\ell )\to
CE_\bullet (g\ell (A),\ell )\; . \]

This map is neither associative nor commutative.  To show that it
induces an associative and commutative product on the space of
coinvariants we need the following simple observation  Let $A$ be a unital
algebra (not necessarily associative) and let $tr:M_n(A)\to A$ be the standard
``trace''.  Then $\alpha\in [M_n(K), M_n(A)]$ (commutator subspace) iff
$Tr\;\alpha =0$.  Now, since in the above map $Tr(\alpha\oplus\beta
)=Tr(\beta\oplus\alpha )$ and $Tr(\alpha\oplus (\beta\oplus\gamma )) =
Tr((\alpha\oplus\beta )\oplus\gamma )$, we conclude that the induced
product
\[ CE_\bullet (g\ell (A,m))_{g\ell (K)}\otimes CE_\bullet (g\ell
(A,m))_{g\ell (K)}\to CE_\bullet (g\ell (A,m))_{g\ell (K)} \]
is graded commutative and associative and hence $CE_\bullet (g\ell
(A),\ell )$ is a DG cocommutative and commutative Hopf algebra.

The rest of the proof of the following extension of Loday-Quillen-Tsygan
theorem is very similar to the original case of associative algebras [L].
We note that an extension of Loday-Quillen-Tsygan theorem to the
category of $DG$ Lie algebras has already been obtained by Burghelea
[B].

\begin{theorem}
Let $(A,m)$ be a unital $A_\infty$-algebra over a field of characteristic
zero.  Then there is a canonical isomorphism of graded Hopf algebras
\[ H_\bullet (g\ell (A),\ell)\simeq\Lambda (HC_\bullet (A,m)[1])\; . \]
\end{theorem}

\vskip30pt

\noindent Masoud Khalkhali\\
University of Western Ontario\\
London, Canada\\
N6A 5B7\\
masoud@julian.uwo.ca

\end{document}